\documentclass[12pt, reqno]{amsart}
\allowdisplaybreaks
\usepackage{bm,amssymb,amsthm,amsmath}
\usepackage[numbers,sort&compress]{natbib}
\newcommand{\upcite}[1]{\textsuperscript{\textsuperscript{\cite{#1}}}}
\title[]{Global well-posedness of 2D nonlinear Boussinesq equations with mixed partial viscosity and thermal diffusivity}
\author{Chao Chen,\,Jitao Liu}
\address[Chao Chen]{School of Mathematics and Computer Science, Fujian Normal University,
Fuzhou, Fujian 350108, P. R. China.}
\email{chenchao\_math@sina.cn, chenchao@fjnu.edu.cn}
\address[Jitao Liu]{College of Applied Sciences, Beijing University of Technology, Beijing, 100124, P. R. China.}
\email{jtliu@bjut.edu.cn,\,\,\,jtliumath@qq.com}

\keywords{nonlinear Boussinesq equations; global well-posedness; partial viscosity; partial thermal diffusivity.}
\thanks{{\em 2010 Mathematics Subject Classification.} 35B30, 35D05, 35Q30, 76D05, 86A04.}

\theoremstyle{plain}
\newtheorem{corollary}{Corollary}[section]
\newtheorem{theorem}{Theorem}[section]
\newtheorem{lemma}{Lemma}[section]
\newtheorem{proposition}{Proposition}[section]

\theoremstyle{definition}
\newtheorem{definition}{Definition}[section]

\newtheorem{remark}{Remark}[section]

\let\f=\frac
\let\p=\partial
\def\R{\Bbb R}

\def\no{\noindent}
\def\endproof{\hphantom{MM}\hfill\llap{$\square$}\goodbreak}
\newcommand{\beq}{\begin{equation}}
\newcommand{\eeq}{\end{equation}}
\newcommand{\ben}{\begin{eqnarray}}
\newcommand{\een}{\end{eqnarray}}
\newcommand{\beno}{\begin{eqnarray*}}
\newcommand{\eeno}{\end{eqnarray*}}

\begin{document}


\begin{abstract}
In this paper, we discuss with the global well-posedness of 2D anisotropic nonlinear Boussinesq equations with any two positive
viscosities and one positive thermal diffusivity. More precisely, for three kinds of viscous combinations, we obtain the global well-posedness without any assumption on the solution. For other three difficult cases, under the minimal regularity assumption, we also derive the unique global solution. To the authors' knowledge, our result is new even for the simplified model, that is, $F(\theta)=\theta e_2$.
\end{abstract}
\maketitle

\section{Introduction and main results}\hspace*{\parindent}
The two dimensional anisotropic nonlinear Boussinesq system in the whole space can be read as:
\begin{equation}\label{cauchy}
\left\{\begin{array}{ll}
\p_t{u^x}+{u}\cdot\nabla {u^x}-\nu_{xx}\p_{xx}{u^x}-\nu_{xy}\p_{yy}{u^x}+\p_x\pi={F_1}(\theta),\\
\p_t{u^y}+{u}\cdot\nabla {u^y}-\nu_{yx}\p_{xx}{u^y}-\nu_{yy}\p_{yy}{u^y}+\p_y\pi={F_2}(\theta),\\
\p_t\theta+{u}\cdot\nabla \theta-\kappa_x\p_{xx}\theta-\kappa_y\p_{yy}\theta=0,\\
\nabla\cdot {u}=0,\\
({u},\theta)(x,0)=({u}_0,\theta_0)(x),
\end{array}\right.
\end{equation}

where ${u}=(u^x,u^y)$, $\theta$ represent the velocity field and temperature respectively. $\pi$ is a scalar pressure, and $\nu_{ij},\,\kappa_i,\,(i,j=x,y)$ are viscosities and thermal diffusivities. The term ${F}(\theta)=(F_1(\theta), F_2(\theta))\in C^2$ is a vector field function satisfying ${F}(0)=0$. For convenience, we use the notation
$$
(A \mid B)=\left( {\begin{array}{*{20}c}
 \nu_{xx} &\nu_{xy}&\vline &\kappa_x   \\
\nu_{yx} &\nu_{yy}&\vline &\kappa_y \\
\end{array} } \right),
$$
as a matrix denoting all the viscosities and thermal diffusivities.

The Boussinesq system arises from the description of natural convection, modeling many geophysical flows such as atmospheric fronts and ocean circulations
(see for example \cite{Majda2,Ped}). Mathematically, it provides an accurate approximation to the 3D incompressible fluids in many important applications.  For example, the 2D Boussinesq equations retain some key features of the 3D Euler and Navier-Stokes equations such as the vortex stretching effects which was shown in \cite{Majda}. It is well known that whether classical solutions to the 3D Euler and Navier-Stokes equations can develop finite time singularities is still open. Therefore, the understanding of 2D Boussinesq equations may shed light on this challenging problem.

Now, let us recall the development of global well-posedness to (\ref{cauchy}) briefly when $F(\theta)=\theta e_2$ with $e_2=(0,1)$. The global regularity of (\ref{cauchy}) with full viscosity and thermal diffusivity is established in \cite{candib}. If all the parameters are zero, (\ref{cauchy}) becomes inviscid and the global regularity problem appears to be out of reach.

Nevertheless, in the past years, the intermediate cases when some of the parameters are positive have attracted considerable attention and big progress has been made. In the beginning, when $\nu_{i,j}\equiv\nu>0$ and $\kappa_x=\kappa_y=0$,  Chae \cite{Chae1} and Hou-Li \cite{HL} establish the global regularity which resolves one of the open problems proposed by Moffatt in \cite{moffatt}. In addition, the case of $\nu_{i,j}\equiv0,\,\kappa_x=\kappa_y>0$ is also studied in \cite{Chae1}. Later, Danchin and Paicu \cite{Danchin} derive the global well-posedness for the anisotropic Boussinesq equations with only horizontal viscosity or thermal diffusivity, i.e. $\nu_{xx}=\nu_{yx}>0,\,\nu_{xy}=\nu_{yy}=\kappa_x=\kappa_y=0$ or $\kappa_x>0,\,\nu_{i,j}=\kappa_y=0$. Recently, Cao and Wu obtain the global regularity for the system with vertical viscosity and thermal diffusivity ($\nu_{xx}=\nu_{yx}=\kappa_x=0,\,\nu_{xy}=\nu_{yy}=\kappa_y>0$ ) in \cite{CaoWu3}. Afterwards, their result is extended to more general source term $F(\theta)$ instead of $\theta e_2$ by Wu and Zheng \cite{WZ}. In a word, they get the unique global solution to the Cauchy problem of (\ref{cauchy}) with
 $$
\hbox{Case 1.\,\,}(A \mid B)=\left( {\begin{array}{*{20}c}
 1 &0&\vline &1   \\
 1 &0&\vline &0 \\
\end{array} } \right)\quad\hbox{and}\quad
\hbox{Case 2.\,\,}(A \mid B)=\left( {\begin{array}{*{20}c}
 0 &1&\vline &0   \\
 0 &1&\vline &1 \\
\end{array} } \right).
$$
Subsequently, Adhikari et al. in \cite{ACSWXY} proved that, when $\nu_{xy}=\nu_{yx}>0,\,\nu_{xx}=\nu_{yy}=\kappa_x=\kappa_y=0$ or $\nu_{yx}=\nu_{yy}>0,\,\nu_{xx}=\nu_{xy}=\kappa_x=\kappa_y=0$, the corresponding system with $F(\theta)=\theta e_2$ always possess global classical solutions. Based on this result, it is not hard to conclude the global well-posedness to the Cauchy problem of (\ref{cauchy}) with
 $$
\hbox{Case 3.\,\,}(A \mid B)=\left( {\begin{array}{*{20}c}
 0 &0&\vline &0   \\
 1 &1&\vline &1 \\
\end{array} } \right),\quad\quad\quad
\hbox{Case 4.\,\,}(A \mid B)=\left( {\begin{array}{*{20}c}
 0 &0&\vline &1   \\
 1 &1&\vline &0 \\
\end{array} } \right),
$$

 $$
\hbox{Case 5.\,\,}(A \mid B)=\left( {\begin{array}{*{20}c}
 0 &1&\vline &1   \\
 1 &0&\vline &0 \\
\end{array} } \right),\quad\quad\quad
\hbox{Case 6.\,\,}(A \mid B)=\left( {\begin{array}{*{20}c}
 0 &1&\vline &0   \\
 1 &0&\vline &1 \\
\end{array} } \right).
$$

Recently, Du and Zhou in \cite{DZ} study the global well-posedness to the MHD system with some kinds of mixed partial viscosities in the whole space. Motivated by this work, in the present paper, we intend to study the global well-posedness for all the rest cases to (\ref{cauchy}) with partial viscosities. It is emphasized that this work covers the global well-posedness of (\ref{cauchy}) with any two positive
viscosities and one positive thermal diffusivity, even for the simplified case $F(\theta)=\theta e_2$. To be more precise, we will study (\ref{cauchy}) with all the following cases,
 $$
\hbox{Case 7.\,\,}(A \mid B)=\left( {\begin{array}{*{20}c}
 1 &1&\vline &1   \\
 0 &0&\vline &0 \\
\end{array} } \right),\quad\quad\quad
\hbox{Case 8.\,\,}(A \mid B)=\left( {\begin{array}{*{20}c}
 0 &1&\vline &1   \\
 0 &1&\vline &0 \\
\end{array} } \right),
$$

 $$
\hbox{Case 9.\,\,}(A \mid B)=\left( {\begin{array}{*{20}c}
 1 &0&\vline &0   \\
 1 &0&\vline &1 \\
\end{array} } \right),\quad\quad\quad
\hbox{Case 10.\,\,}(A \mid B)=\left( {\begin{array}{*{20}c}
 1 &1&\vline &0   \\
 0 &0&\vline &1 \\
\end{array} } \right),
$$

$$
\hbox{Case 11.\,\,}(A \mid B)=\left( {\begin{array}{*{20}c}
 1 &0&\vline &1   \\
 0 &1&\vline &0 \\
\end{array} } \right),\quad\quad\quad
\hbox{Case 12.\,\,}(A \mid B)=\left( {\begin{array}{*{20}c}
 1 &0&\vline &0   \\
 0 &1&\vline &1 \\
\end{array} } \right).
$$

Our main results can be summarized by the following two theorems.

\begin{theorem}\label{gs}
Suppose that $\theta_0\in L^\infty(\R^2)\cap H^2(\R^2)$, $F(\theta_0)\in L^2(\R^2)$ and $u_0\in H^2(\R^2)$ with $\nabla\cdot u_0=0$. Then there exists a unique global solution $(u,\theta)$ solving the system (\ref{cauchy}) with Case 7--Case 9 such that
\beno
(u,\theta)\in L^\infty(\R_+;H^{2}(\R^2)).
\eeno
\end{theorem}
\begin{remark}\label{dif0}
Due to lack of partial viscosities and thermal diffusivities, it makes this problem more challenging.
Mathematically, without sufficient smooth effect, $H^1$ and $H^2$ estimates cannot be obtained by the standard energy methods. To overcome this difficulty, we make full advantage of the incompressible condition together with the anisotropic inequalities to obtain some delicate estimates, which helps us to get the global well-posedness. On the other hand, there are six separate cases addressed in this theorem. To make the proof more clearly and avoid repetition, we try the best to use the uniform proof.
\end{remark}

\begin{theorem}\label{gs2}
Let $\theta_0\in L^\infty(\R^2)\cap H^2(\R^2)$, $F(\theta_0)\in L^2(\R^2)$ and $u_0\in H^2(\R^2)$ with $\nabla\cdot u_0=0$.

(i)\,Then there exists a unique global solution $(u,\theta)$ solving the system (\ref{cauchy}) with Case 10--Case 11 such that $(u,\theta)\in L^\infty(0,T;H^{2}(\R^2))$ for any $T>0$, if
\ben\label{con1}
\|\p_xu^y\|_{L^2(0,T;L^2(\R^2))}<\infty\quad\quad\hbox{or}\quad\quad\|\p_x\theta\|_{L^2(0,T;L^2(\R^2))}<\infty.
\een

(ii)\,Then there exists a unique global solution $(u,\theta)$ solving the system (\ref{cauchy}) with Case 12 such that $(u,\theta)\in L^\infty(0,T;H^{2}(\R^2))$ for any $T>0$, if
\ben\label{con2}
\|\p_yu^x\|_{L^2(0,T;L^2(\R^2))}<\infty\quad\quad\hbox{or}\quad\quad\|\p_y\theta\|_{L^2(0,T;L^2(\R^2))}<\infty.
\een
\end{theorem}
\begin{remark}\label{dif}
Compared with other cases, Case 11 and 12 are more difficult to prove since less smooth effect is provided. Essentially, the a priori estimates of the vorticity cannot be obtained without additional assumptions.
Though there are two viscosities, only $\|\p_xu^x\|_{L^2(0,T;L^2(\R^2))}$ can be obtained by energy methods. Furthermore, more challenge occurs in the $H^2$ estimates of velocity, since the viscosity term is not enough to control all the nonlinear terms. To obtain the global well-posedness
under minimal regularity assumptions, we need more delicate estimates.
\end{remark}
\begin{remark}
The solutions in Theorem \ref{gs} and \ref{gs2} are defined as follows.
\end{remark}
\begin{definition}\label{weak}
Suppose $\theta_0\in L^\infty\cap H^1(\R^2)$,\,\,$F(\theta_0)\in L^2(\R^2)$ and $u_0\in H^1(\R^2)$. A pair of measurable $u(x,t)$ and $\theta(x,t)$  is called a global weak solution of (\ref{cauchy}) if for any $T>0$,
\ben\label{weak1}
&&u\in L^\infty(0,T;H^{2}(\R^2)),\,\,\theta\in L^\infty(0,T;H^{2}(\R^2));
\een
and
\ben\label{weak2}
&&\int_{\R^2}u_0\cdot\varphi_0\,dx+\int_0^T\int_{\R^2}\big[u\cdot\varphi_t+u\cdot\nabla \varphi\cdot u+F(\theta)\cdot\varphi\big]dxdydt\\
&=&\sum_{i,j\in\{x,y\}}\nu_{ij}\int_0^T\int_{\R^2}\partial_ju^i\partial_j\varphi^idxdydt,\notag
\een
\ben\label{weak3}
\int_{\Omega}\theta_0\psi_0dx+\int_0^T\int_{\R^2}\big[\theta\psi_t+u\cdot\nabla\psi\theta\big]dxdydt=\sum_{i\in\{x,y\}}\kappa_{i}\int_0^T\int_{\R^2}\partial_i\theta\partial_i\psi dxdydt,
\een
holds for any $\varphi=(\varphi^x,\varphi^y),\,\psi\in C^\infty([0,T]\times\R^2)$ with $\nabla\cdot\varphi=\varphi(x,T)=\psi(x,T)=0$.
\end{definition}

\begin{remark}
The existence of such weak solution is based on the a priori estimates  and the Friedrichs approximation method. See Lemma \ref{gw} for details.
\end{remark}

This paper is organized as follows. In Section 2, we introduce some notations and technical lemmas used for our estimates in the rest sections. In Section 3, we  establish the a priori estimates. Section 4 is devoted to the global well-posedness (i.e. the proof of Theorem \ref{gs} and \ref{gs2}).

\section{Preliminary}\hspace*{\parindent}

In this section, we introduce some definitions and useful lemmas throughout this paper. To begin with it, we define the inner products on $L^2(\R^2)$ and space $V$ by
$$(u,v)=\sum\limits_{{i}\in\{x,y\}}\int_{\R^2}u^iv^i\,dxdy,$$
and
$$V=\{u\in H^1(\R^2): \nabla\cdot u = 0\,\,\hbox{in}\,\,\R^2\},$$
respectively.  We also denote the dual space of $V$ by $V'$  and the action of $V'$ on $V$ by $<\cdot\,,\,\cdot>$. Moreover, we set the trilinear continuous form by
\ben\label{b0}
b(u,v,w)=\sum\limits_{{i,j}\in\{x,y\}}\int_{\R^2}u^i\p_iv^jw^j\,dx.
\een
If $u\in\,V$, it is obvious that
\ben\label{b1}
b(u,v,w)=-b(u,w,v),\,\,\,\forall\,\,v,w\in H^1(\R^2),
\een
and
\ben\label{b2}
b(u,v,v)=0,\,\,\,\forall\,\,v\in H^1(\R^2).
\een
\begin{remark}\label{gw1}
One can follow standard arguments as in the theory of the Navier-Stokes equations (see e.g., \cite{RT}) to conclude that the system (\ref{weak2})-(\ref{weak3}) is equivalent to the following system
\begin{equation}\label{weak4}
\f{d}{dt}<u,\varphi>+\sum_{i,j\in\{x,y\}}\nu_{ij}(\partial_ju^i,\partial_j\varphi^i)=(F(\theta),\varphi)-b(u,u,\varphi),
\end{equation}

\begin{equation}\label{weak5}
\f{d}{dt}<\theta,\psi>+\sum_{i\in\{x,y\}}\kappa_{i}(\partial_i\theta,\partial_i\psi)=-b(u,\theta,\psi),
\end{equation}
for any $\varphi\in L^2(0,T; V)$ and $\psi\in L^2(0,T; H^1(\R^2))$.
\end{remark}

Thanks to the incompressible condition on $u$, we can obtain the following uniform bound on $\theta$ with the space variables.
\begin{lemma}\label{lem0}
Assume that $(u,\,\theta)$ is a smooth solution of system $(\ref{cauchy})$ and $\theta_0\in L^p(\R^2)$ with $p\in [2,\infty]$, then there holds
\beno
&&\|\theta(t,\cdot)\|_{L^p(\R^2)}\leq\ C\|\theta_0\|_{L^p(\R^2)},
\eeno
for any $t>0$, where $C$ is an absolute constant.
\end{lemma}

Since $F(\theta)$ is a $C^2$ function on the variable $\theta$, we can conclude that
\ben\label{M}
\|{F(\theta(t,\cdot))}\|_{L^{\infty}(\R^2)}+\|{F'(\theta(t,\cdot))}\|_{L^{\infty}(\R^2)}+\|{F''(\theta(t,\cdot))}\|_{L^{\infty}(\R^2)}\leq\ M,
\een
from Lemma \ref{lem0}. Here $M$ is a constant depending on $\|\theta_0\|_{L^\infty}$ only.

Now, let us recall some classical results which can be found in the cited references.
\begin{lemma}\label{lem:vor}\upcite{Majda}
 Suppose that the vector field $u$ is divergent free. The vorticity $w=\p_yu^x-\p_xu^y\in L^p(\R^2)\,\,\hbox{for}\,\,p\in(1,\infty)$, then there exists a constant $C_p$ depending on $p$  such that
 \beno
 \|\nabla u\|_{L^p(\R^2)}\leq C_p\|w\|_{L^p(\R^2)}.
 \eeno
 \end{lemma}

\begin{lemma}\label{elliptic}\upcite{GT}
Consider the following elliptic system
\beno
-\Delta f=g\quad&\hbox{in}\,\,\R^2.
\eeno
If $g\in L^{p}(\R^2),$ then for any $p\in(1,\infty)$, there exists a unique solution $f$ satisfying the estimates
$$\|\nabla^2f\|_{L^{p}(\R^2)}\leq C_p\|g\|_{L^{p}(\R^2)},$$
where $C$ depending only on $p$.
\end{lemma}

As an application of Lemma \ref{elliptic}, $-\Delta u^x=\p_yw$ and $-\Delta u^y=\p_xw$, the following lemma holds.

\begin{lemma}\label{lem:curl}
 Suppose that $u=(u^x,u^y)$ is an vector field, whose divergence is zero. Its vorticity $w=\p_yu^x-\p_xu^y\in L^p(\R^2)\,\,\hbox{for}\,\,p\in(1,\infty)$, then there exists a constant $C_p$ only depending on $p$ such that
 \beno
\|\nabla^2 u^x\|_{L^p(\R^2)}+\|\p_y\nabla u\|_{L^p(\R^2)}\leq C_p\|\p_yw\|_{L^p(\R^2)},
 \eeno
  \beno
\|\nabla^2 u^y\|_{L^p(\R^2)}+\|\p_x\nabla u\|_{L^p(\R^2)}\leq C_p\|\p_xw\|_{L^p(\R^2)},
 \eeno
  \beno
\|\nabla^2 \p_yu^x\|_{L^p(\R^2)}+\|\p_{yy}\nabla u\|_{L^p(\R^2)}\leq C_p\|\p_{yy}w\|_{L^p(\R^2)},
 \eeno
  \beno
\|\nabla^2 \p_xu^y\|_{L^p(\R^2)}+\|\p_{xx}\nabla u\|_{L^p(\R^2)}\leq C_p\|\p_{xx}w\|_{L^p(\R^2)},
 \eeno
   \beno
\|\nabla^2 \p_xu^x\|_{L^p(\R^2)}\leq C_p\|\p_{xy}w\|_{L^p(\R^2)}.
 \eeno
 \end{lemma}

Then, we have to mention here that the following two classical inequalities also play the key role in the proof of Theorem \ref{gs} and \ref{gs2}.

\begin{lemma}\label{lem:p2}[Lemma 2.2, \cite{CaoWu3}]
Assume $f,g,h$ are smooth functions in $\R^2$. Then it holds that
 \beno
\int_{\R^2}|fgh|\,dxdy\leq
C\|f\|_{L^2(\R^2)}\|g\|_{L^2(\R^2)}^{\f12}\|\p_xg\|_{L^2(\R^2)}^{\f12}\|h\|_{L^2(\R^2)}^{\f12}\|\p_yh\|_{L^2(\R^2)}^{\f12},
\eeno
where  $C$ is an absolute constant.
\end{lemma}

\begin{lemma}\label{lem:s2}[Lemma A.2, \cite{Danchin}]
For any smooth function $f(x,y)$, there exists a constant C such that
 \beno
\|f\|_{L^\infty(\R^2)}\leq C(\|f\|_{L^2(\R^2)}+\|\p_xf\|_{L^2(\R^2)}+\|\p_{yy}f\|_{L^2(\R^2)}),
 \eeno
and
\beno
\|f\|_{L^\infty(\R^2)}\leq C(\|f\|_{L^2(\R^2)}+\|\p_yf\|_{L^2(\R^2)}+\|\p_{xx}f\|_{L^2(\R^2)}).
\eeno
\end{lemma}

Before the proof of our theorems, we would like to point out the following basic facts.
\begin{proposition}\label{symmetry}
Any smooth solution $(u,\theta)$ of (\ref{cauchy}) also solves following system
\begin{equation}\label{symmetry1}
\left\{\begin{array}{ll}
\p_t{u^x}+{u}\cdot\nabla {u^x}-\nu_{yy}\p_{xx}{u^x}-\nu_{yx}\p_{yy}{u^x}+\p_x\pi={F_2}(\theta),\\
\p_t{u^y}+{u}\cdot\nabla {u^y}-\nu_{xy}\p_{xx}{u^y}-\nu_{xx}\p_{yy}{u^y}+\p_y\pi={F_1}(\theta),\\
\p_t\theta+{u}\cdot\nabla \theta-\kappa_y\p_{xx}\theta-\kappa_x\p_{yy}\theta=0,\\
\nabla\cdot {u}=0,\\
({u},\theta)(x,0)=({u}_0,\theta_0)(x),
\end{array}\right.
\end{equation}
\end{proposition}
\begin{remark}\label{symmetry2}
Setting $(U^y,U^x,\Theta,\Pi)(x,y,t)=(u^x,u^y,\theta,\pi)(y,x,t)$ and putting it into (\ref{cauchy}) imply Proposition \ref{symmetry} directly.
\end{remark}
\begin{remark}\label{symmetry3}
Since $F_1,F_2$ have the same regularity, it is clear that to prove Theorem \ref{gs} and \ref{gs2}, it suffices to deal with Case 7,\,8,\,10 and 11.
\end{remark}

\section{A priori Estimates}\hspace*{\parindent}
In this section, we will establish the a priori estimates.
\begin{lemma}\label{lem1}
Suppose that $\theta_0\in L^2\cap L^\infty(\R^2)$, $F(\theta_0)\in L^2(\R^2)$ and $u_0\in L^2(\R^2)$ with $\nabla\cdot u_0=0$. Then for a smooth solution $(u,\theta)$ of (\ref{cauchy}) with $\nu_{i,j},\,\kappa_i\geq 0\,(i,j=x,y)$, the following estimates hold
\ben\label{energy1}
&&\|\theta\|_{L^2}^2+\int_0^T({\kappa_x}\|\p_x \theta\|_{L^2}^2+{\kappa_y}\|\p_y \theta\|_{L^2}^2)dt\leq\ C,
\een
\ben\label{energy2}
\|{F}(\theta)\|_{L^2}^2+\int_0^T({\kappa_x}\|\p_x{F}(\theta)\|_{L^2}^2+{\kappa_y}\|\p_y{F}(\theta)\|_{L^2}^2)dt\leq\ C(T),
\een
\ben\label{energy3}
\|{u}\|_{L^2}^2+\int_0^T[({\nu_{xx}}+{\nu_{yy}})\|\p_xu^x\|_{L^2}^2+{\nu_{xy}}\|\p_yu^x\|_{L^2}^2+{\nu_{yx}}\|\p_xu^y\|_{L^2}^2]dt\leq\ C(T),
\een
where $C$ is an absolute constant and the constant $C(T)$ depends only on T.
\end{lemma}
\no{\bf Proof.}\quad
The proof of (\ref{energy1}) relies on standard methods and we omit it here. To get (\ref{energy2}), by multiplying ${F'}(\theta)$ on both sides of $(\ref{cauchy})^{3}$, we obtain that ${F}(\theta)$ solves the following equation
\ben\label{eqF}
\p_t{F}(\theta)+{u}\cdot\nabla{F}(\theta)-\kappa_x\p_{xx}{F}(\theta)-\kappa_y\p_{yy}{F}(\theta)=-{F''}(\theta)[\kappa_x(\p_x\theta)^2+\kappa_y(\p_y\theta)^2].
\een

Taking inner product of (\ref{eqF}) with ${F}(\theta)$ and integrating on $\R^2$ yields that
\beno
&&\f12\f{d}{dt}\|{F}(\theta)\|_{L^2}^2+{\kappa_x}\|\p_x{F}(\theta)\|_{L^2}^2+{\kappa_y}\|\p_y{F}(\theta)\|_{L^2}^2\\
&\leq&C\|{F''}(\theta)\|_{L^\infty}[{\kappa_x}\|\p_x\theta\|_{L^2}^2+{\kappa_y}\|\p_y\theta\|_{L^2}^2]\\
&\leq&CM[{\kappa_x}\|\p_x\theta\|_{L^2}^2+{\kappa_y}\|\p_y\theta\|_{L^2}^2],
\eeno
This, together with integrating in time and (\ref{energy1}), implies that
\ben\label{21}
\|{F}(\theta)\|_{L^2}^2+\int_0^T({\kappa_x}\|\p_x{F}(\theta)\|_{L^2}^2+{\kappa_y}\|\p_y{F}(\theta)\|_{L^2}^2)dt\leq\ C(T).
\een

Through similar process, one can also obtain
\beno
&&\f12\f{d}{dt}\|{u}\|_{L^2}^2+[({\nu_{xx}}+{\nu_{yy}})\|\p_xu^x\|_{L^2}^2+{\nu_{xy}}\|\p_yu^x\|_{L^2}^2+{\nu_{yx}}\|\p_xu^y\|_{L^2}^2]\\
&\leq&C\|{F}(\theta)\|_{L^2}\|u\|_{L^2}^2,
\eeno
which yields the conclusion by turning to the Gronwall's inequality and (\ref{21}).
\endproof

\begin{proposition}\label{lem3}
Under the assumptions of Theorem \ref{gw}, for a smooth solution $(u,\theta)$ of (\ref{cauchy}) with Case 7, there holds
\ben\label{estimate31}
\|{u}\|_{H^2}^2+\int_0^T[\|\p_{xy}w\|_{L^2}^2+\|{\p_{yy}}w\|_{L^2}^2]dt\leq\ C(T),
\een
and
\ben\label{estimate32}
\|{\theta}\|_{H^2}^2+\int_0^T\|\nabla^2\p_x\theta\|_{L^2}^2dt\leq\ C(T),
\een
where the constant $C(T)$ depends only on T.
\end{proposition}
\no{\bf Proof.}\quad
By Lemma \ref{lem1}, it holds that
\ben\label{estimate33}
\|\theta\|_{L^2}^2+\|{u}\|_{L^2}^2+\int_0^T[\|\p_x \theta\|_{L^2}^2+\|\nabla u^x\|_{L^2}^2]dt\leq C(T).
\een

{\bf Step 1. $H^1$ estimates}\\

By taking $\p_i\,(i=x,y)$ on both sides of the second equation in $(\ref{cauchy})$, we have
\ben\label{31}
\p_t\p_i\theta+u\cdot\nabla \p_i\theta-\p^2_x\p_i\theta=-\p_iu\cdot\nabla\theta.
\een
Multiplying (\ref{31}) with $\p_i\theta$ and integrating on $\R^2$ imply that
\begin{equation}\label{32}
\f12\f{d}{dt}\|\nabla \theta\|_{L^2}^2+\|\p_x\nabla\theta\|_{L^2}^2=-\sum\limits_{{i,j}\in\{x,y\}}\int_{\R^2}\p_iu^j\p_j\theta\p_i\theta\,dxdy\triangleq\sum\limits_{{i,j}\in\{x,y\}}I^{ij},
\end{equation}
where $I^{ij}=-\int_{\R^2}\p_iu^j\p_j\theta\p_i\theta\,dxdy$.

We directly apply Lemma \ref{lem:p2} and Young inequality to obtain that
\begin{equation}\label{33}
\begin{array}{lll}
|I^{ix}|&\leq C\|\nabla u^x\|_{L^2}\|\p_x\theta\|_{L^2}^{\f12}\|\p_{xy}\theta\|_{L^2}^{\f12}\|\nabla\theta\|_{L^2}^{\f12}\|\p_x\nabla\theta\|_{L^2}^{\f12}\\
&\leq C\|\partial_x\nabla\theta\|_{L^2}\|\nabla u^x\|_{L^2}\|\nabla\theta\|_{L^2}\\
&\leq\dfrac{1}{8}\|\p_x\nabla\theta\|_{L^2}^2+C\|\nabla u^x\|_{L^2}^2\|\nabla\theta\|_{L^2}^2.
\end{array}
\end{equation}
In addition,  noting the incompressible condition on $u$ yields that
\ben\label{34}
|I^{iy}|&=&|\int_{\R^2}\p_xu^y\p_x\theta\p_y\theta\,dxdy|+|\int_{\R^2}\p_xu^x\p_y\theta\p_y\theta\,dxdy|\notag\\
&=&|\int_{\R^2}u^y\p_y\theta\p_{xx}\theta\,dxdy|+|\int_{\R^2}u^y\p_{xy}\theta\p_x\theta\,dxdy|+2|\int_{\R^2}u^x\p_y\theta\p_{xy}\theta\,dxdy|\notag\\
&\leq&C\|\nabla\p_{x}\theta\|_{L^2}\|u\|_{L^2}^{\f12}\|\nabla u^x\|_{L^2}^{\f12}\|\nabla\theta\|_{L^2}^{\f12}\|\nabla\partial_x\theta\|_{L^2}^{\f12}\notag\\
&\leq&\f{1}4\|\p_x\nabla\theta\|_{L^2}^2+C\|u\|_{L^2}^2\|\nabla u^x\|_{L^2}^2\|\nabla\theta\|_{L^2}^2.
\een
It follows from (\ref{estimate33}), (\ref{32})-(\ref{34}) and Gronwall's inequality that
\ben\label{estimate36}
\|\nabla\theta\|_{L^2}^2+\int_0^T\|\p_x\nabla\theta\|_{L^2}^2dt\leq C(T).
\een

To obtain the higher estimates on $u$ , we are ready to do the estimate of vorticity. Taking $\nabla\times$ to the equation of velocity in $(\ref{cauchy})$ gives that
\ben\label{41}
\p_tw+u\cdot\nabla w-\Delta\p_yu^x=\p_x(F_2(\theta))-\p_y(F_1(\theta)),
\een
where $w=\p_xu^y-\p_yu^x$. Then, multiplying (\ref{41}) with $w$ and integrating on $\R^2$ yield that
\beno
&&\f12\f{d}{dt}\|w\|_{L^2}^2+\|\nabla\p_yu\|_{L^2}^2\\
&\leq&\int_{\R^2}F'_2(\theta)\p_x\theta w\,dx-\int_{\R^2}F'_1(\theta)\p_y\theta w\,dx\\
&\leq&C(\|F'_1(\theta)\|_{L^\infty}\|\p_y\theta\|_{L^2}+\|F'_2(\theta)\|_{L^\infty}\|\p_x\theta\|_{L^2})\|w\|_{L^2}\\
&\leq&C\|\nabla\theta\|_{L^2}^2+\|w\|_{L^2}^2.
\eeno
Finally, the Gronwall's inequality guarantees that
\ben\label{estimate37}
\|w\|_{L^2}^2+\int_0^T\|\nabla\p_yu\|_{L^2}^2dt\leq\ C(T).
\een
{\bf Step 2. $H^2$ estimates}\\

Multiplying (\ref{41}) with $-\Delta w$ and integrating on $\R^2$ yield that
\ben\label{w0}
&&\f12\f{d}{dt}\|\nabla w\|_{L^2}^2+\|\p_{xy}w\|_{L^2}^2+\|{\p_{yy}}w\|_{L^2}^2\notag\\
&=&-\int_{\R^2}\nabla w\cdot\nabla u\cdot\nabla w\,dx-\int_{\R^2}\p_x(F_2(\theta))\Delta w\,dx+\int_{\R^2}\p_y(F_1(\theta))\Delta w\,dx\notag\\
&=&\sum_i^3J^i.
\een
One can make use of Lemma \ref{lem:p2} to obtain that
\ben\label{w1}
J^1&\leq&|\int_{\R^2}\p_x w\p_x u^x\p_x w\,dxdy|+|\int_{\R^2}\p_x w\p_x u^y\p_y w\,dxdy|\notag\\
&&+|\int_{\R^2}\p_y w\p_y u^x\p_x w\,dxdy|+|\int_{\R^2}\p_y w\p_y u^y\p_y w\,dxdy|\notag\\
&\leq&C\|\p_xu^x\|_{L^2}^{\f12}\|\p_{xx}u^x\|_{L^2}^{\f12}\|\p_xw\|_{L^2}^{\f32}\|\p_{xy}w\|_{L^2}^{\f12}\\
&&+C\|\p_xu^y\|_{L^2}^{\f12}\|\p_{xy}u^y\|_{L^2}^{\f12}\|\p_xw\|_{L^2}\|\p_yw\|_{L^2}^{\f12}\|\p_{xy}w\|_{L^2}^{\f12}\notag\\
&&+C\|\p_yu^x\|_{L^2}^{\f12}\|\p_{xy}u^x\|_{L^2}^{\f12}\|\p_xw\|_{L^2}^{\f12}\|\p_yw\|_{L^2}\|\p_{xy}w\|_{L^2}^{\f12}\notag\\
&&+C\|\p_yu^y\|_{L^2}^{\f12}\|\p_{yy}u^y\|_{L^2}^{\f12}\|\p_yw\|_{L^2}^{\f32}\|\p_{xy}w\|_{L^2}^{\f12}\notag\\
&\leq&C\|w\|_{L^2}^{\f12}\|\p_{xy}u\|_{L^2}^{\f12}\|\nabla w\|_{L^2}^{\f32}\|\p_{xy}w\|_{L^2}^{\f12}\notag\\
&\leq&\f14\|\p_{xy}w\|_{L^2}^2+C\|w\|_{L^2}^{\f23}\|\p_{xy}u\|_{L^2}^{\f23}\|\nabla w\|_{L^2}^2.\notag
\een

Thanks to integration by parts, Lemma \ref{lem:p2}, H\"{o}lder and Young inequalities, we conclude that the second term $J^2$ satisfies
\ben\label{w2}
J^2&=&\int_{\R^2}{\p_{xx}}(F_2(\theta))\p_xw\,dxdy-\int_{\R^2}\p_y(F_2(\theta))\p_{xy}w\,dxdy\notag\\
&=&\int_{\R^2}F''_2(\theta)\p_x\theta\p_x\theta\p_xw+F'_2(\theta){\p_{xx}}\theta\p_xw\,dxdy-\int_{\R^2}F'_2(\theta)\p_y\theta\p_{xy}w\,dxdy\notag\\
&\leq&C\|\p_xw\|_{L^2}\|\p_x\theta\|_{L^2}\|{\p_{xx}}\theta\|_{L^2}^{\f12}\|\p_{xy}\theta\|_{L^2}^{\f12}+C\|\p_xw\|_{L^2}\|{\p_{xx}}\theta\|_{L^2}\\
&&+C\|\p_y\theta\|_{L^2}\|\p_{xy}w\|_{L^2}\notag\\
&\leq&\f14\|\p_{xy}w\|_{L^2}^2+C(1+\|\p_x\theta\|_{L^2}^2)\|\p_xw\|_{L^2}^2+C(\|\p_y\theta\|_{L^2}^2+\|\nabla\p_x\theta\|_{L^2}^2).\notag
\een

Similarly, the third term $J^3$ satisfies the estimate as follows.
\ben\label{w3}
J^3&=&-\int_{\R^2}\p_{xy}(F_1(\theta))\p_xw\,dxdy+\int_{\R^2}F'_1(\theta)\p_y\theta{\p_{yy}}w\,dxdy\notag\\
&=&-\int_{\R^2}F''_1(\theta)\p_x\theta\p_y\theta\p_xw\,dxdy+\int_{\R^2}F'_1(\theta)\p_{xy}\theta\p_xw\,dxdy\notag\\
&&+\int_{\R^2}F'_1(\theta)\p_y\theta{\p_{yy}}w\,dxdy\notag\\
&\leq&C\|\p_xw\|_{L^2}\|\p_x\theta\|_{L^2}^{\f12}\|\p_y\theta\|_{L^2}^{\f12}\|\p_{xy}\theta\|_{L^2}+C\|\p_xw\|_{L^2}\|\p_{xy}\theta\|_{L^2}\\
&&+C\|\p_y\theta\|_{L^2}\|{\p_{yy}}w\|_{L^2}\notag\\
&\leq&\f14\|{\p_{yy}}w\|_{L^2}^2+C(1+\|\nabla\theta\|_{L^2}^2)\|\p_xw\|_{L^2}^2+C(\|\p_y\theta\|_{L^2}^2+\|\nabla\p_x\theta\|_{L^2}^2).\notag
\een

Hence, it holds that
 \beno
&&\f{d}{dt}\|\nabla w\|_{L^2}^2+\|\p_{xy}w\|_{L^2}^2+\|{\p_{yy}}w\|_{L^2}^2\\
&\leq&C(1+\|\nabla\theta\|_{L^2}^2+\|w\|_{L^2}^{\f23}\|\p_{xy}u\|_{L^2}^{\f23})\|\nabla w\|_{L^2}^2+C(\|\nabla\theta\|_{L^2}^2+\|\nabla\p_x\theta\|_{L^2}^2),
\eeno
This, together with (\ref{estimate36}), (\ref{estimate37}) and Gronwall's inequality, leads to
\ben\label{estimate38}
\|\nabla w\|_{L^2}^2+\int_0^T[\|\p_{xy}w\|_{L^2}^2+\|{\p_{yy}}w\|_{L^2}^2]dt\leq C(T).
\een

The next thing is to give the $H^2$ estimates of $\theta$. Multiplying $(\ref{31})^1$ with $-\p^3_x\theta$ and $(\ref{31})^2$ with $-\p^3_y\theta$ yield that
\ben\label{t0}
&&\f12\f{d}{dt}\|\Delta \theta\|_{L^2}^2+\|\Delta\p_{x}\theta\|_{L^2}^2\notag\\
&=&\int_{\R^2}u\cdot\nabla\p_x\theta\p^3_x\theta\,dx+\int_{\R^2}u\cdot\nabla\p_y\theta\p^3_y\theta\,dx\\
&&-\int_{\R^2}\p_xu\cdot\nabla\theta\p^3_x\theta\,dx-\int_{\R^2}\p_yu\cdot\nabla\theta\p^3_y\theta\,dx\notag\\
&=&\sum_i^4H^i\notag
\een
We conclude that
\ben\label{t1}
&&H^3+H^1\notag\\
&=&-\int_{\R^2}\p_xu^x\p_{x}\theta\p^3_x\theta\,dx-\int_{\R^2}\p_xu^y\p_{y}\theta\p^3_x\theta\,dx\notag\\
&&+\int_{\R^2}u^x{\p_{xx}}\theta\p^3_x\theta\,dx+\int_{\R^2}u^y\p_{xy}\theta\p^3_x\theta\,dx\notag\\
&\leq&C\|\p^3_x\theta\|_{L^2}\|\p_xu^x\|_{L^2}^{\f12}\|\p_{xy}u^x\|_{L^2}^{\f12}\|\p_{x}\theta\|_{L^2}^{\f12}\|\p^2_{x}\theta\|_{L^2}^{\f12}\\
&&+C\|\p^3_x\theta\|_{L^2}\|\p_xu^y\|_{L^2}^{\f12}\|\p_{xy}u^y\|_{L^2}^{\f12}\|\p_{y}\theta\|_{L^2}^{\f12}\|\p_{xy}\theta\|_{L^2}^{\f12}\notag\\
&&+C\|\p^3_x\theta\|_{L^2}\|u^x\|_{L^2}^{\f12}\|\p_yu^x\|_{L^2}^{\f12}\|\p_{xx}\theta\|_{L^2}^{\f12}\|\p^3_{x}\theta\|_{L^2}^{\f12}\notag\\
&&+C\|\p^3_x\theta\|_{L^2}\|u^y\|_{L^2}^{\f12}\|\p_xu^y\|_{L^2}^{\f12}\|\p_{xy}\theta\|_{L^2}^{\f12}\|\p_{xyy}\theta\|_{L^2}^{\f12}\notag\\
&\leq&C\|\Delta\p_x\theta\|_{L^2}\|w\|_{L^2}^{\f12}\|\p_{xy}u\|_{L^2}^{\f12}\|\nabla\theta\|_{L^2}^{\f12}\|\Delta\theta\|_{L^2}^{\f12}\notag\\
&&+C\|\Delta\p_x\theta\|_{L^2}^{\f32}\|u\|_{L^2}^{\f12}\|w\|_{L^2}^{\f12}\|\Delta\theta\|_{L^2}^{\f12}\notag\\
&\leq&\f14\|\Delta\p_x\theta\|_{L^2}^2+C(\|\nabla\theta\|_{L^2}^2+\|u\|_{L^2}^2\|w\|_{L^2}^2)\|\Delta \theta\|_{L^2}^2\notag\\
&&+C\|w\|_{L^2}^2\|\p_{xy}u\|_{L^2}^2,\notag
\een
by H\"{o}lder inequalities, Lemma \ref{lem:p2} and Young inequalities.

As for the left terms, more careful decompositions are needed. Thanks to Lemma \ref{lem:p2}, it holds that
\ben\label{t2}
&&H^2\notag\\
&=&\int_{\R^2}u^x\p_{xy}\theta\p^3_y\theta\,dx+\int_{\R^2}u^y{\p_{yy}}\theta\p^3_y\theta\,dx\notag\\
&=&-\int_{\R^2}\p_yu^x\p_{xy}\theta\p_{yy}\theta\,dx-\int_{\R^2}u^x\p_{xyy}\theta\p_{yy}\theta\,dx-\f12\int_{\R^2}\p_yu^y\p_{yy}\theta\p_{yy}\theta\,dx\notag\\
&\leq&C\|\p_{xy}\theta\|_{L^2}^{\f12}\|\p_{xyy}\theta\|_{L^2}^{\f12}\|\p_yu^x\|_{L^2}\|\p_{yy}\theta\|_{L^2}^{\f12}\|\p_{xyy}\theta\|_{L^2}^{\f12}\\
&&+C\|\p_{xyy}\theta\|_{L^2}\|u^x\|_{L^2}^{\f12}\|\p_yu^x\|_{L^2}^{\f12}\|\p_{yy}\theta\|_{L^2}^{\f12}\|\p_{xyy}\theta\|_{L^2}^{\f12}\notag\\
&&+C\|\p_{yy}\theta\|_{L^2}^{\f32}\|\p_{xyy}\theta\|_{L^2}^{\f12}\|\p_yu^y\|_{L^2}^{\f12}\|\p_{yy}u^y\|_{L^2}^{\f12}\notag\\
&\leq&C\|\Delta\p_{x}\theta\|_{L^2}\|w\|_{L^2}\|\Delta\theta\|_{L^2}+C\|\Delta\p_{x}\theta\|_{L^2}^{\f32}\|u\|_{L^2}^{\f12}\|w\|_{L^2}^{\f12}\|\Delta\theta\|_{L^2}^{\f12}\notag\\
&&+C\|\Delta\p_{x}\theta\|_{L^2}^{\f12}\|\Delta\theta\|_{L^2}^{\f32}\|w\|_{L^2}^{\f12}\|\p_{xy}u\|_{L^2}^{\f12}\notag\\
&\leq&\f14\|\Delta\p_x\theta\|_{L^2}^2+C\|w\|_{L^2}^2(1+\|u\|_{L^2}^2)\|\Delta\theta\|_{L^2}^2+C\|w\|_{L^2}^{\f23}\|\p_{xy} u\|_{L^2}^{\f23}\|\Delta\theta\|_{L^2}^2,\notag
\een
and
\ben\label{t3}
H^4&=&-\int_{\R^2}\p_yu^x\p_{x}\theta\p^3_y\theta\,dx-\int_{\R^2}\p_yu^y\p_{y}\theta\p^3_y\theta\,dx\notag\\
&=&\int_{\R^2}\p_yu^x\p_{xy}\theta\p_{yy}\theta\,dx+\int_{\R^2}\p_{yy}u^x\p_{x}\theta\p_{yy}\theta\,dx\notag\\
&&+\int_{\R^2}\p_{yy}u^y\p_{y}\theta\p_{yy}\theta\,dx+\int_{\R^2}\p_yu^y\p_{yy}\theta\p_{yy}\theta\,dx\notag\\
&\leq&C\|\p_{xy}\theta\|_{L^2}^{\f12}\|\p_{yy}\theta\|_{L^2}^{\f12}\|\p_{xyy}\theta\|_{L^2}\|w\|_{L^2}\notag\\
&&+C\|\p_{yy}u^x\|_{L^2}\|\p_{yy}\theta\|_{L^2}^{\f12}\|\p_{xyy}\theta\|_{L^2}^{\f12}\|\p_{x}\theta\|_{L^2}^{\f12}\|\p_{xy}\theta\|_{L^2}^{\f12}\\
&&+C\|\p_{yy}u^y\|_{L^2}\|\p_{yy}\theta\|_{L^2}^{\f12}\|\p_{xyy}\theta\|_{L^2}^{\f12}\|\p_y\theta\|_{L^2}^{\f12}\|\p_{yy}\theta\|_{L^2}^{\f12}\notag\\
&&+C\|\p_{yy}\theta\|_{L^2}\|\p_yu^y\|_{L^2}^{\f12}\|{\p_{yy}}u^y\|_{L^2}^{\f12}\|\p_{yy}\theta\|_{L^2}^{\f12}\|\p_{xyy}\theta\|_{L^2}^{\f12}\notag\\
&\leq&C\|\Delta\p_{x}\theta\|_{L^2}\|\Delta\theta\|_{L^2}\|w\|_{L^2}+C\|\Delta\p_{x}\theta\|_{L^2}^{\f12}\|\nabla w\|_{L^2}\|\nabla\theta\|_{L^2}^{\f12}\|\Delta\theta\|_{L^2}\notag\\
&&+C\|\Delta\p_{x}\theta\|_{L^2}^{\f12}\|\Delta\theta\|_{L^2}^{\f32}\|w\|_{L^2}^{\f12}\|\nabla w\|_{L^2}^{\f12}\notag\\
&\leq&\f14\|\Delta\p_x\theta\|_{L^2}^2+C(1+\|w\|_{L^2}^2+\|\nabla w\|_{L^2}^2+\|w\|_{L^2}^2\|\nabla w\|_{L^2}^2)\|\Delta\theta\|_{L^2}^2\notag\\
&&+C\|\nabla\theta\|_{L^2}^2.\notag
\een
Now summing up (\ref{t0})-(\ref{t3}), employing Gronwall's inequality together with (\ref{estimate36}), (\ref{estimate37}) and (\ref{estimate38}) give that
\ben\label{estimate39}
\|\Delta \theta\|_{L^2}^2+\int_0^T\|\Delta\p_{x}\theta\|_{L^2}^2dt\leq C(T),
\een
which leads to (\ref{estimate31})-(\ref{estimate32}) by combing (\ref{estimate36}), (\ref{estimate37}), (\ref{estimate38})and applying Lemma \ref{elliptic}, \ref{lem:curl}.
\endproof

For Case 8,\,10, we have the similar estimates.
\begin{corollary}\label{cor5}
Assume that the assumptions of Theorem \ref{gw} hold.\\
(i)\,For  a smooth solution $(u,\theta)$ of (\ref{cauchy}) with Case 8, it holds that
\ben\label{estimate51}
\|{u}\|_{H^2}^2+\int_0^T[\|\p_{xy}w\|_{L^2}^2+\|{\p_{yy}}w\|_{L^2}^2]dt\leq\ C(T),
\een
and
\ben\label{estimate52}
\|{\theta}\|_{H^2}^2+\int_0^T\|\nabla^2\p_x\theta\|_{L^2}^2dt\leq\ C(T).
\een
(ii)\,For a smooth solution $(u,\theta)$ of (\ref{cauchy}) with Case 10, if (\ref{con1}) holds, then we have
\ben\label{estimate91}
\|{u}\|_{H^2}^2+\int_0^T[\|\p_{xy}w\|_{L^2}^2+\|{\p_{yy}}w\|_{L^2}^2]dt\leq\ C(T),
\een
and
\ben\label{estimate92}
\|{\theta}\|_{H^2}^2+\int_0^T\|\nabla^2\p_y\theta\|_{L^2}^2dt\leq\ C(T),
\een
where the constant $C(T)$ only depends on T.
\end{corollary}
\no{\bf Proof.}\quad
(i)\,\,\,Compared with the proof of Lemma \ref{lem3}, there is only a minor modification in the $H^1$ estimates of $\theta$.
Different from (\ref{estimate33}) , by Lemma \ref{lem1}, it holds here that
\ben\label{estimate53}
\|{u}\|_{L^2}^2+\int_0^T\|\p_y u\|_{L^2}^2dt\leq\ C(T).
\een
Moreover, by same calculation with (\ref{31}), (\ref{32}) and the incompressible condition, we can obtain that
\ben\label{51}
&&\f12\f{d}{dt}\|\nabla \theta\|_{L^2}^2+\|\p_x\nabla\theta\|_{L^2}^2\notag\\
&\leq&\f{1}2\|\p_x\nabla\theta\|_{L^2}^2+C\|\p_y u\|_{L^2}^2\|\nabla\theta\|_{L^2}^2+C\|u\|_{L^2}^2\|\p_y u\|_{L^2}^2\|\nabla\theta\|_{L^2}^2.
\een
Then one can follow Lemma \ref{lem3} to finish the proof.\\
(iii)\,\,\,Case 10 is similar to Case 7 except horizontal thermal diffusivity replaced by vertical thermal diffusivity.

Since the symmetric structure of (\ref{cauchy}), it suffices to show the $H^1$ estimates of $\theta$.  It is clear that the energy estimates provide
\ben\label{estimate93}
\|{u}\|_{L^2}^2+\int_0^T\|\nabla u^x\|_{L^2}^2\leq\ C(T).
\een
Compared with (\ref{32}),  it is only different in the first term. Since (\ref{con1}) holds, it follows that
\ben\label{91}
I^1&=&2|\int_{\R^2}u^y\p_x\theta\p_{xy}\theta\,dx|\notag\\
&\leq&C\|u^y\|_{L^2}\|\p_x u^y\|_{L^2}^{\f12}\|\p_x\theta\|_{L^2}^{\f12}\|\p_{xy}\theta\|_{L^2}^{\f32}\\
&\leq&\f{1}4\|\p_{xy}\theta\|_{L^2}^2+C\|u^y\|_{L^2}^2\|\p_x u^y\|_{L^2}^2\|\p_x\theta\|_{L^2}^2.\notag
\een
Then, the $H^1$ estimates of $\theta$ can be obtained in a similar way.
\endproof

As mentioned in Remark \ref{dif}, although there are two viscosities in Case 11, only one viscosity function takes effect by the incompressible condition. It brings less viscosity terms than other cases in the a priori estimates. Therefore, we need an additional assumption to overcome this difficulty.
\begin{corollary}\label{cor10}
Suppose that the assumptions of Theorem \ref{gw} hold, for a smooth solution $(u,\theta)$ of (\ref{cauchy}) with Case 11, the estimates
\ben\label{estimate101}
\|{u}\|_{H^2}^2+\|{\theta}\|_{H^2}^2+\int_0^T[\|\p_{xy}w\|_{L^2}^2+\|\nabla^2\p_x\theta\|_{L^2}^2]dt\leq\ C(T),
\een
hold if (\ref{con1}) be satisfied, where the constant $C(T)$ only depends on T.
\end{corollary}
\no{\bf Proof.}\quad
Assume that (\ref{con1}) holds, by the similar calculation with Proposition \ref{lem3}, it is clear to obtain that
\ben\label{estimate102}
\|w\|_{L^2}^2+\|\nabla\theta\|_{L^2}^2+\int_0^T[\|\p_{xy}u\|_{L^2}^2+\|\p_x\nabla\theta\|_{L^2}^2]dt\leq C(T).
\een
As for the $H^2$ estimates, we have less useful information than the other cases. It seems hard to obtain the single estimate of $\|w\|_{L^\infty(0,T;L^2(\R^3))}$ or $\|\Delta\theta\|_{L^\infty(0,T;L^2(\R^3))}$ separately  as before. Hence, we intend to estimate them together.
It follows from integration by parts that
\ben\label{103}
&&\f12\f{d}{dt}(\|\nabla w\|_{L^2}^2+\|\Delta \theta\|_{L^2}^2)+\|\p_{xy}w\|_{L^2}^2+\|\Delta\p_{x}\theta\|_{L^2}^2\notag\\
&=&\sum_i^3J^i+\sum_i^4H^i.
\een
It should be noted that $J^1+J^2$ and $H^1+H^2+H^3$ can be estimated as before. Therefore, we have
\ben\label{104}
&&J^1+J^2\notag\\
&\leq&\f14\|\p_{xy}w\|_{L^2}^2+C(1+\|\p_x\theta\|_{L^2}^2+\|w\|_{L^2}^{\f23}\|\p_{xy}u\|_{L^2}^{\f23})\|\nabla w\|_{L^2}^2\\
&&+C(\|\p_y\theta\|_{L^2}^2+\|\nabla\p_x\theta\|_{L^2}^2),\notag
\een
and
\ben\label{105}
&&H^1+H^2+H^3\notag\\
&\leq&\f14\|\Delta\p_x\theta\|_{L^2}^2+C(\|w\|_{L^2}^2+\|u\|_{L^2}^2\|w\|_{L^2}^2+\|w\|_{L^2}^{\f23}\|\p_{xy} u\|_{L^2}^{\f23})\|\Delta \theta\|_{L^2}^2\\
&&+C(1+\|\nabla\theta\|_{L^2}^2)\|\Delta \theta\|_{L^2}^2+C\|w\|_{L^2}^2\|\p_{xy}u\|_{L^2}^2.\notag
\een

On the other hand, $J^3$ and $H^4$ can be controlled by the bounds in the following inequalities
\ben\label{106}
J^3&=&\int_{\R^2}\p_{x}(F_1(\theta))\p_{xy}w\,dx-\int_{\R^2}{\p_{yy}}(F_1(\theta))\p_{y}w\,dx\notag\\
&=&\int_{\R^2}F'_1(\theta)\p_x\theta\p_{xy}w\,dx-\int_{\R^2}F''_1(\theta)\p_y\theta\p_y\theta\p_yw\,dx\notag\\
&&-\int_{\R^2}F'_1(\theta)\p_{yy}\theta\p_yw\,dx\notag\\
&\leq&C\|\p_x\theta\|_{L^2}\|\p_{xy}w\|_{L^2}+C\|\p_yw\|_{L^2}\|\p_y\theta\|_{L^2}\|\p_{xy}\theta\|_{L^2}^{\f12}\|\p_{yy}\theta\|_{L^2}^{\f12}\\
&&+C\|\p_yw\|_{L^2}\|\p_{yy}\theta\|_{L^2}\notag\\
&\leq&\f14[\|\p_{xy}w\|_{L^2}^2+\|\Delta\theta\|_{L^2}^2]+C(1+\|\nabla\theta\|_{L^2}^2)(1+\|\nabla w\|_{L^2}^2),\notag
\een
and
\ben\label{107}
H^4&=&-\int_{\R^2}\p_yu^x\p_{x}\theta\p^3_y\theta\,dx-\int_{\R^2}\p_yu^y\p_{y}\theta\p^3_y\theta\,dx\notag\\
&=&\int_{\R^2}\p_yu^x\p_{xy}\theta\p_{yy}\theta\,dx+\int_{\R^2}\p_{yy}u^x\p_{x}\theta\p_{yy}\theta\,dx\notag\\
&&+\int_{\R^2}\p_yu^y\p_{yy}\theta\p_{yy}\theta\,dx+\int_{\R^2}\p_{yy}u^y\p_{y}\theta\p_{yy}\theta\,dx\notag\\
&\leq&C\|\p_{xy}\theta\|_{L^2}^{\f12}\|\p_{yy}\theta\|_{L^2}^{\f12}\|\p_{xyy}\theta\|_{L^2}\|w\|_{L^2}\notag\\
&&+C\|\p_{yy}u^x\|_{L^2}\|\p_{yy}\theta\|_{L^2}^{\f12}\|\p_{xyy}\theta\|_{L^2}^{\f12}\|\p_{x}\theta\|_{L^2}^{\f12}\|\p_{xy}\theta\|_{L^2}^{\f12}\notag\\
&&+C\|\p_{yy}\theta\|_{L^2}\|\p_yu^y\|_{L^2}^{\f12}\|{\p_{yy}}u^y\|_{L^2}^{\f12}\|\p_{yy}\theta\|_{L^2}^{\f12}\|\p_{xyy}\theta\|_{L^2}^{\f12}\\
&&+C\|\p_{yy}u^y\|_{L^2}\|\p_y\theta\|_{L^2}^{\f12}\|\p_{yy}\theta\|_{L^2}^{\f12}\|\p_{yy}\theta\|_{L^2}^{\f12}\|\p_{xyy}\theta\|_{L^2}^{\f12}\notag\\
&\leq&C\|\Delta\p_{x}\theta\|_{L^2}\|\Delta\theta\|_{L^2}\|w\|_{L^2}+C\|\nabla w\|_{L^2}\|\Delta\theta\|_{L^2}^{\f12}\|\Delta\p_{x}\theta\|_{L^2}^{\f12}\|\nabla\theta\|_{L^2}^{\f12}\|\nabla\p_x\theta\|_{L^2}^{\f12}\notag\\
&&+C\|\Delta\p_{x}\theta\|_{L^2}^{\f12}\|\Delta\theta\|_{L^2}(\|\Delta\theta\|_{L^2}^{\f12}\|w\|_{L^2}^{\f12}\|\p_{xy} u\|_{L^2}^{\f12}+\|\p_{xy} u\|_{L^2}\|\nabla\theta\|_{L^2}^{\f12})\notag\\
&\leq&\f14\|\Delta\p_x\theta\|_{L^2}^2+C(1+\|w\|_{L^2}^2+\|\p_{xy} u\|_{L^2}^2+\|w\|_{L^2}^2\|\p_{xy} u\|_{L^2}^2)\|\Delta\theta\|_{L^2}^2\notag\\
&&+C\|\nabla\theta\|_{L^2}^2(1+\|\Delta\theta\|_{L^2}^2)+C\|\nabla\p_x\theta\|_{L^2}\|\nabla w\|_{L^2}^2.\notag
\een

Now, the estimates are closed by turning to the Gronwall's inequality and the proof is finished.
\endproof

\section{Global well-posedness}\hspace*{\parindent}

After obtaining the necessary estimates, we can prove the existence and uniqueness of the solution in Theorem \ref{gs} and Theorem \ref{gs2}.
\begin{proposition}\label{gw}
Let $\theta_0\in L^\infty(\R^2)\cap H^2(\R^2)$,\,\,$F(\theta_0)\in L^2(\R^2)$ and $u_0\in H^2(\R^2)$ with $\nabla\cdot u_0=0$. Then there exists at least one global weak solution $(u,\theta)$ solving  system (\ref{cauchy}) with Case 7, Case 8, Case 10 and Case 11 such that
\beno
(u,\theta)\in L^\infty(\R_+;H^{2}(\R^2)).
\eeno
\end{proposition}
\no{\bf Proof.}\quad
Based on our a priori estimates, the proof can be achieved through Friedrichs method which is also known as ``modified Galerkin method". One may see \cite{Taylor} for example. Without loss of generality, we assume $\nu_{ij}=\kappa_i=1\,(i,j=x,y)$.
For $\epsilon>0$, let $j$ be a positive radial compactly supported smooth function whose integral equals 1 and denote $J_\epsilon$ as a Friedrichs mollifier by
\beno
J_\epsilon=j_\epsilon\ast u,\quad\,\hbox{where}\,\,\,\,\,j_\epsilon=\epsilon^{-2}j(\epsilon^{-1}x).
\eeno
Moreover, let $\mathcal{P}$ be the Leray projector over divergence free vector fields. Then the following properties holds
\ben\label{jp}
J^2_\epsilon=J_\epsilon,\quad\,\mathcal{P}^2=\mathcal{P},\quad\,\mathcal{P}J_\epsilon=J_\epsilon\mathcal{P}.
\een

Now, consider the following approximating system
\begin{equation}\label{cauchy2}
\left\{\begin{array}{ll}
\p_tu_\epsilon+\mathcal{P}J_\epsilon(J_\epsilon u_\epsilon\cdot\nabla J_\epsilon u_\epsilon)-\mathcal{P}J_\epsilon (F(J_\epsilon\theta_\epsilon))
=\Delta\mathcal{P}J_\epsilon u_\epsilon,\\
\p_t\theta_\epsilon+J_\epsilon(J_\epsilon u_\epsilon\cdot\nabla J_\epsilon\theta_\epsilon)=\Delta J_\epsilon\theta_\epsilon,\\
(u_\epsilon,\theta_\epsilon)(x,0)=J_\epsilon(u_0,\theta_0).
\end{array}\right.
\end{equation}
The Cauchy-Lipschitz theorem guarantees the existence of a unique smooth solution $(u_\epsilon,\,\theta_\epsilon)$ in short time. Due to (\ref{jp}),  $(\mathcal{P}u_\epsilon,\,\theta_\epsilon)$ and $(J_\epsilon u_\epsilon,\,J_\epsilon\theta_\epsilon)$ are also solutions of (\ref{cauchy2}). This yields that the solution to \eqref{cauchy2} exactly solves the following system
\begin{equation}\label{cauchy3}
\left\{\begin{array}{ll}
\p_tu_\epsilon+\mathcal{P}J_\epsilon (u_\epsilon\cdot\nabla u_\epsilon)-\mathcal{P}J_\epsilon (F(\theta_\epsilon))
=\Delta\mathcal{P}u_\epsilon,\\
\p_t\theta_\epsilon+J_\epsilon(u_\epsilon\cdot\nabla\theta_\epsilon)=\Delta \theta_\epsilon,\\
\nabla\cdot u_\epsilon=0,\\
(u,\theta)(x,0)=J_\epsilon(u_0,\theta_0).
\end{array}\right.
\end{equation}

Then, by Lemmas \ref{lem:vor}-\ref{lem:curl}, the fact that $J_\epsilon$ and $\mathcal{P}J_\epsilon$ are orthogonal projectors in $L^2$ and the similar priori estimates obtained in Section 3, we have
\ben\label{epsilon}
(u_\epsilon,\theta_\epsilon)\in L^\infty(\R_+;H^{2}(\R^2)),
\een
where the uniform bound here is independent of $\epsilon$.

Thanks to the Sobolev embeddings and H\"{o}lder inequalities, $u_\epsilon\theta_\epsilon\in L_{\rm{loc}}^2(\R_+;L^{4}(\R^2))$ and $u_\epsilon\otimes u_\epsilon\in L_{\rm{loc}}^2(\R_+;L^{2}(\R^2))$. These imply that $(\p_tu_\epsilon\,\,\p_t\theta_\epsilon)\in L_{\rm{loc}}^2(\R_+;H^{-1}(\R^2)).$ Since the Sobolev embeddings $L^4\hookrightarrow H^{-1}$ and $L^2\hookrightarrow H^{-1}$ are locally compact, the classical Aubin-Lions argument (with a diagonal process if needed) guarantees that we can extract a sequence of $(u_\epsilon,\,\,\theta_\epsilon)$, still denoted by itself, such that there is a limit $(u,\theta)$ satisfying $u^\epsilon\longrightarrow u\,\,\hbox{strongly\,\,in}\,\,L^2$ and $\theta^\epsilon\longrightarrow \theta\,\,\hbox{strongly\,\,in}\,\,L^2$. Moreover,
\beno
(u,\theta)\in L^\infty(\R_+;H^{2}(\R^2)).
\eeno
Hence, the existence is obtained. The left thing is to show the uniqueness of the solution.

For any fixed $T>0$, suppose there are two solutions $(u,\theta,\pi)$, $(\widetilde{u},\widetilde{\theta},\widetilde{\pi})$ of (\ref{cauchy}) and let $U=\widetilde{u}-u,\,\,\Theta=\widetilde{\theta}-\theta,\,\,\Pi=\widetilde{\pi}-\pi,$ then by Remark \ref{gw1}, it holds that
\ben\label{unix}
&&\f{d}{dt}<U,\varphi>+b(\widetilde{u},U,\varphi)+b(U,u,\varphi)+\sum_{i,j\in\{x,y\}}\nu_{ij}(\partial_jU^i,\partial_j\varphi^i)\notag\\
&&=([F(\widetilde{\theta})-F(\theta)],\varphi),
\een
\ben\label{uniy}
\f{d}{dt}<\Theta,\psi>+b(\widetilde{u},\theta,\psi)+b(U,\theta,\varphi)+\sum_{i\in\{x,y\}}\kappa_{i}(\partial_i\Theta,\partial_i\psi)=0,
\een
\ben\label{uniz}
(U,B)(x,0)=0,
\een
for any $\varphi\in L^2(0,T; V)$ and $\psi\in L^2(0,T; H^1(\R^2))$.

Now, we take $\varphi=U$ in $(\ref{unix})$ and $\psi=\Theta$ in $(\ref{uniy})$ respectively. By (\ref{b2}) and Lions-Magenes Lemma (see e.g., \cite{RT}), one can obtain that
\ben\label{uni0}
&&\f12\f{d}{dt}(\|U\|_{L^2}^2+\|\Theta\|_{L^2}^2)\\
&\leq&\int_{\R^2}[F(\widetilde{\theta})-F(\theta)]\cdot U\,dx-\int_{\R^2}U\cdot\nabla u\cdot U\,dx-\int_{\R^2}U\cdot\nabla \theta\cdot \Theta\,dx\notag
\een

From (\ref{M}), it is clear that
\ben\label{uni00}
\|F(\widetilde{\theta})-F(\theta)\|_{L^2}\leq\|\int_{\theta}^{\widetilde{\theta}}F'(s)\,ds\|_{L^2}\leq C\|\Theta\|_{L^2}.
\een
Thus, we have
\beno
&&\f{d}{dt}(\|U\|_{L^2}^2+\|\Theta\|_{L^2}^2)\\
&\leq&C(1+\|\nabla u\|_{L^\infty}^2+\|\nabla \theta\|_{L^\infty}^2)(\|U\|_{L^2}^2+\|\Theta\|_{L^2}^2).
\eeno
This, together with Gronwall's inequality, implies that
\ben\label{uni000}
\|U(t)\|_{L^2}^2+\|\Theta(t)\|_{L^2}^2\leq Ce^{\int_0^T[\|\nabla u\|_{L^\infty}+\|\nabla \theta\|_{L^\infty}]dt}(\|U_0\|_{L^2}^2+\|\Theta_0\|_{L^2}^2).
\een

Next, we recall the a priori estimates obtained in Section 3. By applying Lemma \ref{lem:curl} and \ref{lem:s2}, it holds that
\beno
\|\nabla u\|_{L^\infty}
&\leq&2\|\p_x u^x\|_{L^\infty}+\|\p_x u^y\|_{L^\infty}+\|\p_y u^x\|_{L^\infty}\notag\\
&\leq&C(\|\p_x u^x\|_{L^2}+\|\p_{xx} u^x\|_{L^2}+\|\p_{xyy}u^x\|_{L^2})\\
&&+C(\|\p_x u^y\|_{L^2}+\|\p_{xx} u^y\|_{L^2}+\|\p_{xyy}u^y\|_{L^2})\notag\\
&&+C(\|\p_y u^x\|_{L^2}+\|\p_{yy} u^x\|_{L^2}+\|\p_{xxy}u^x\|_{L^2})\notag\\
&\leq&C(\|w\|_{L^2}+\|\nabla w\|_{L^2}+\|\p_{xy}w\|_{L^2}),\notag
\eeno
and
\beno
\|\nabla \theta\|_{L^\infty}
\leq C(\|\nabla\theta\|_{L^2}+\|\Delta \theta\|_{L^2}+\|\Delta\p_{x}\theta\|_{L^2}),\notag
\eeno
or
\beno
\|\nabla \theta\|_{L^\infty}
\leq C(\|\nabla\theta\|_{L^2}+\|\Delta \theta\|_{L^2}+\|\Delta\p_{y}\theta\|_{L^2}).\notag
\eeno

Then it follows Proposition \ref{lem3}, Corollary \ref{cor5}, Corollary \ref{cor10} and (\ref{uni000}) that
\beno
\|U(t)\|_{L^2}^2+\|\Theta(t)\|_{L^2}^2\leq C(T)(\|U_0\|_{L^2}^2+\|\Theta_0\|_{L^2}^2)=0,
\eeno
for any $t\in[0,T].$  Hence, it is done and the proof of Theorem \ref{gs} and \ref{gs2} is complete.
\endproof

\section*{Acknowledgments}
Chen is supported by NSFC 11301079 and Fujian educational program.
\par

\end{document}